\documentclass{gtpart}
\usepackage{amsmath,amssymb}
\usepackage[dvips]{graphics}
\usepackage{amsfonts,amscd,latexsym,bbm,epsfig,epic,eepic,oldgerm,psfrag}

\title{Symplectic Manifolds with Vanishing Action--Maslov Homomorphism}

\author{Mark Branson}
\givenname{Mark}
\surname{Branson}
\address{Department of Mathematics\\
The Technion, Israeli Institute of Technology\\
Haifa, 32000\\
Israel}
\email{branson@tx.technion.ac.il}

\arxivreference{0910.5259}
\arxivpassword{rgdkz}

\volumenumber{}
\issuenumber{}
\publicationyear{}
\papernumber{}
\startpage{}
\endpage{}
\doi{}
\MR{}
\Zbl{}
\received{}
\revised{}
\accepted{}
\published{}
\publishedonline{}
\proposed{}
\seconded{}
\corresponding{}
\editor{}
\version{}

\newcommand{\univ}{{\rm univ}}

\newcommand{\id}{{\rm id}}
\newcommand{\CP}{{\mathbb CP}}
\newcommand{\Ham}{{\rm Ham}}

\newcommand{\ev}{{\rm ev}}

\newtheorem{theorem}{Theorem}[section]

\newtheorem{corollary}[theorem]{Corollary}

\newtheorem{lemma}[theorem]{Lemma}
\newtheorem{proposition}[theorem]{Proposition}

\newtheorem{note}[theorem]{Note}
\newtheorem{definition}[theorem]{Definition}

\numberwithin{figure}{section}
\numberwithin{equation}{section}
\numberwithin{table}{section}

\newcommand{\QH}{{\rm QH}}
\newcommand{\eff}{{\rm eff}}
\newcommand{\enr}{{\rm enr}}
\newcommand{\proofend}{\hfill$\square$}

\begin{document}
\begin{abstract}
The action--Maslov homomorphism $I\co\pi_1(\Ham(X,\omega))\to\R$ is an important tool for understanding the topology of the Hamiltonian group of monotone symplectic manifolds.  We explore conditions for the vanishing of this homomorphism, and show that it is identically zero when the Seidel element has finite order and the homology satisfies property $\mathcal{D}$ (a generalization of having homology generated by divisor classes).  We use these results to show that $I=0$ for products of projective spaces and the Grassmannian of $2$ planes in $\C^4$.  
\end{abstract}
\maketitle

\section{Introduction}\label{introduction}
Let $(X, \omega)$ be a monotone symplectic manifold, that is $\omega = \kappa c_1$ on $\pi_2(X)$ with $\kappa > 0$. Polterovich introduced the action--Maslov homomorphism $I \co  \pi_1 (\Ham(X, \omega)) \to \R$ in \cite{polterovich}. Manifolds where $I = 0$ have many interesting properties. When $I=0$, the Hamiltonian group has infinite Hofer diameter \cite{polterovich}, the asymptotic spectral invariants descend from $\widetilde{\text{Ham}}(X,\omega)$ to $\text{Ham}(X,\omega)$ \cite{mcduff10}, the Calabi quasimorphism descends from $\widetilde{\text{Ham}}(X,\omega)$ to $\text{Ham}(X,\omega)$  \cite{ep09}, and another non-Calabi quasimorphism exists on $\text{Ham}(X,\omega)$ \cite{py}.  Another very important result states that a K\"{a}hler manifold cannot be K\"{a}hler--Einstein unless $I$ vanishes on all holomorphic Hamiltonian circle actions \cite{ep09}, \cite{shelukhin}.

There are a number of manifolds where $I$ is known to vanish:
\begin{itemize}
\item $\CP^2$ \cite{polterovich}
\item $S^2\times S^2$ \cite{polterovich}
\item $\CP^n$ \cite{ep09}, \cite{mcduff10}
\item $2$-manifolds \cite{shelukhin}
\item $\CP^2\# 3\overline{\CP^2}$ \cite{branson}
\end{itemize}

In \cite{mcduff10}, McDuff gives several conditions which imply $I = 0$ (or an equivalent condition in the nonmonotone case). Essentially, these criteria specify manifolds where most of the genus zero Gromov--Witten invariants vanish or have carefully controlled properties. We extend these results by exploring the form of the Seidel element more deeply. For monotone symplectic manifolds, the Seidel element always has integral coefficients and a finite number of terms. By studying these constraints on the Seidel element and properties of the quantum homology, we can show that $I$ vanishes for products of projective spaces and the Grassmannian $G(2, 4)$. 

\begin{theorem}\label{cpnproducts}
 $I = 0$ for $\CP^{n_1} \times \ldots \times \CP^{n_k}$ with a monotone symplectic form. 
 \end{theorem}
 
 \begin{theorem}\label{g24}
 $I = 0$ for $G(2, 4)$ with the monotone symplectic form.
 \end{theorem}

Theorem \ref{cpnproducts} is related to results of Pedroza \cite{pedroza} and Leclercq \cite{leclercq}. They showed that, for $X '$ and $X ''$ monotone symplectic manifolds, $\gamma ' \in \pi_1 (\Ham(X ' ))$, $\gamma '' \in \pi_1 (\Ham(X '' ))$, then $\mathcal{S}(\gamma ' \times \gamma '' ) = \mathcal{S}(\gamma ' ) \otimes \mathcal{S}(\gamma '' )$, where $\mathcal{S}(\gamma)$ is the Seidel element. This is sufficient to show that $I = 0$ for any loop $\gamma ' \in \pi_1 (\Ham(\CP^m \times \CP^n ))$ which is a product of loops in the hamiltonian groups of $\CP^m$ and $\CP^n$. Our result shows that $I = 0$ for \emph{all} loops in $\pi_1 (\Ham(\CP^m \times\CP^n ))$. 

Our method for proving that the action--Maslov homomorphism vanishes depends on showing two facts. First we show that, when $X$ is one of the above manifolds and $\gamma \in \pi_1 (\Ham(X,\omega))$, there exists $k > 0$ such that the Seidel element $\mathcal{S}(k\gamma) = \mathbbm{1} \otimes \lambda$, where $\lambda$ is in the Novikov coefficient ring $\Lambda$ and $\mathbbm{1}$ is the fundamental class $[X]$. $[X]$ is the unit for both the intersection product on $H_*(X)$ and the quantum product on $\QH_*(X)$, and we will use the notation $\mathbbm{1}$ in both. Then, we must prove that $\nu(\mathbbm{1} \otimes \lambda) = 0$, where $\nu$ is the valuation map on quantum homology. These terms will be defined in Section \ref{definitions}, but are well-known in symplectic topology.  The condition we use which is not well-known is Property $\mathcal{D}$ from \cite{mcduff10}.  We say that the quantum homology has property $\mathcal{D}$ if the nontrivial genus zero Gromov-Witten invariants with two homology constraints vanish unless both terms are in the subgroup generated by the divisors or both are in its additive complement.  This condition will also be stated more explicitly in Section \ref{definitions}.

\begin{proposition} \label{nuiszero}
Let $(X, \omega)$ be a symplectic manifold. Suppose that the quantum homology $\QH_* (X, \Lambda)$ has property $\mathcal{D}$. Then for all $\gamma \in \pi_1 (\Ham(X, \omega))$ such that $\mathcal{S}(\gamma) = \mathbbm{1} \otimes \lambda$, $\nu(\mathbbm{1} \otimes \lambda) = 0$. 
\end{proposition}

The following theorem is an immediate consequence.

\begin{theorem}\label{ivanishes}
Let $(X, \omega)$ be a monotone symplectic manifold. Assume that $(X, \omega)$ has property $\mathcal{D}$ and that for all $\gamma \in \pi_1 (\Ham(X, \omega))$, $\exists n$ such that the Seidel element $\mathcal{S}(n\gamma) = \mathbbm{1} \otimes \lambda$ for some $ \lambda \in \Lambda$. Then $I = 0$.
\end{theorem}

These conditions are rather restrictive, but they are satisfied for almost all manifolds where $I$ is known to vanish (with the possible exception of $\CP^2\# 3\overline{\CP^2}$, which may not have $\mathcal{S}(n\gamma) = \mathbbm{1} \otimes \lambda$). 

Property $\mathcal{D}$ is trivially satisfied when the even homology classes are generated by divisors, so it includes many well-studied examples, such as toric varieties.  In many of these cases, it is difficult to show that that $\mathcal{S}(n\gamma) = \mathbbm{1} \otimes \lambda$.  We will say that such Seidel elements have \emph{finite order} (this is not strictly true, but the reason why this is a good term will be discussed in section \ref{definitions}), and that Seidel elements without this property have \emph{infinite order}.  

\noindent\textit{Acknowledgments} I would like to thank my thesis advisor, Dusa McDuff, for inspiring this work and guiding me through the process of writing my thesis and thus this article. Also, I would like to thank Mike Chance, Martin Pinsonnault, Egor Shelukhin, Yaron Ostrover, Leonid Polterovich, and Aleksey Zinger for conversations and correspondence regarding the material herein.  I would also like to thank Michael Entov and the Technion for sponsoring me while I completed this work.  My work was partially supported by the Israel Science Foundation grant No. 723/10.


\section{Definitions}\label{definitions}
Let $X$ be a $2N$ dimensional symplectic manifold.  Most of these definitions can be found in greater detail in \cite{ms}. Let $K^{\eff}$, the effective cone of $(X, \omega)$, be the additive cone generated by the spherical homology classes $A \in H_2^S (X)$ with nonvanishing genus zero Gromov--Witten invariants $\langle a, b, c \rangle^X_A\neq 0$. Consider the Novikov ring $\Lambda_\enr$ given by formal sums 
\begin{equation}\label{novikovring}
\lambda= \sum_{A\in K^\eff} \lambda(A)e^{-A}
\end{equation}

\noindent with the finiteness condition that, $\forall c \in \R$, 
\begin{equation}\label{finiteness}
\#\lbrace A \in K^\eff |\lambda(A) \neq 0, \omega(A) \leq c\rbrace < \infty
\end{equation}

\noindent where $\lambda(A) \in \R$. $\Lambda_\enr$ has a grading given by $|e^A| =2c_1 (A)$. We will call this the enriched Novikov ring. The universal ring $\Lambda_{\text{univ}}$ is $\Lambda [q, q^{-1} ]$, where $\Lambda$ is generated by formal power series of the form: 
\begin{equation}\label{lambdauniv}
\lambda= \sum_{\epsilon\in\R} \lambda_\epsilon t^{-\epsilon}
\end{equation}

\noindent with a similar finiteness condition and $\lambda_\epsilon \in \R$.  The grading on $\Lambda_{\text{univ}}$ is given by setting $\deg(q) = 2$. Note that $\Lambda$ is a field. The map $\varphi \co  \Lambda_\enr \to \Lambda_{\text{univ}}$ is given by taking 
\begin{equation*}
\varphi(e^{-A} ) = q^{-c_1 (A)} t^{-\omega(A)}
\end{equation*}

\noindent and extending by linearity. 

While these two rings are thus related, different properties of the quantum homology become apparent when different coefficient rings are used. In section \ref{cpnsection} the enriched Novikov ring $\Lambda_\enr$ will be used to show that the units are of a specific form for $X$ a product of projective spaces, because the proof fails with the universal ring. In section \ref{nuvanishes}, calculations will be carried out using the universal ring $\Lambda$. The quantum homology with respect to the Novikov ring $\Lambda_\enr$ is given by $\QH_* (X, \Lambda_\enr ) = H_* (X, \R) \otimes \Lambda_\enr $. The grading on $\QH_* (X, \Lambda_\enr )$ will be given by the sum of the grading on $H_* (X, \R)$ and the grading on $\Lambda_\enr $. 

The quantum homology admits a product structure, called the quantum product. Let $\xi_i$ be a basis of $H_* (X)$ and $\xi_i^*$ a dual basis with respect to the intersection product.  Then the quantum product of $a, b \in H_* (X, \R)$ is defined by:
\begin{equation}\label{quantumproduct}
a*b= \sum_{i,A\in K^\eff} \langle a, b, \xi_i\rangle_A^X \xi_i^*\otimes e^{-A}
\end{equation}

We can then extend this to $\QH_* (X, \Lambda_\enr )$ by linearity. The quantum homology $\QH_* (X, \Lambda_{\text{univ}})$ is defined analogously and the map $\id \otimes \varphi$ extends to a ring homomorphism $\Phi \co  \QH_* (X, \Lambda_\enr ) \to \QH_* (X, \Lambda_{\text{univ}})$. We define the valuation map $\nu \co  \QH_* (X, \Lambda_{\text{univ}}) \to \R$ by 
\begin{equation*}
\nu(\sum_{i\in\R} \lambda_i \otimes q^{a_i} t^{b_i} ) = \max\lbrace b_i |\lambda_i \neq 0\rbrace.
\end{equation*}

Next, we will discuss $\mathcal{S}(\gamma)$, the Seidel element (defined in \cite{seidel}). Given a loop $\gamma \in \pi_1 (\Ham(X, \omega))$ with $\gamma = \lbrace\phi_t \rbrace$, we define the Hamiltonian fiber bundle $P_\gamma$ over $S^2$ with fiber $X$. This bundle is given by the clutching construction - take two copies of $D^2 , D_+$ and $D_- $. Then take $X \times D_+$ and glue it to $X \times D_-$ (where $D_-$ has the opposite orientation from $D_+$ ) via the map
\begin{equation*}
(\phi_t (x), e^{2\pi i t} )_+ \sim (x, e^{2\pi it} )_-.
\end{equation*}

When the loop $\gamma$ is clear from context, we will refer to this bundle as $P$. $P_\gamma$ has two canonical classes - the vertical Chern class, denoted $c_1^{\text{vert}}$, and the coupling class, denoted $u_\gamma$. $c_1^{\text{vert}}$ is the first Chern class of the vertical tangent bundle. $u_\gamma$ is the unique class such that $u_\gamma |_X = \omega$ and $u_\gamma^{n+1} = 0$. Given a section class $\sigma \in H_2 (P_\gamma , \Z)$, we can define the Seidel element in $\QH_* (X, \Lambda_\enr )$ by taking 
\begin{equation}\label{seidelenriched}
\mathcal{S}(\gamma, \sigma) = \sum_{A\in H_2^S (X),i} \langle \xi_i \rangle^{P_\gamma}_{\sigma+A} \xi_i^* \otimes e^{-A}
\end{equation}

\noindent where  $H^S_2(X)$ is the image of $\pi_2 (X)$ in $H_2 (X)$ (the spherical homology classes). Note that $\sigma + A$ is a slight abuse of notation; we should actually write $\sigma + \iota_* (A)$, where $\iota \co  X \to P_\gamma$ is the inclusion map. We will continue this abuse throughout the paper. The Seidel element can also be defined in $\QH_* (X, \Lambda_\univ)$ - in this case, the dependence on $\sigma$ is eliminated by an averaging process. 
\begin{equation}\label{seideluniv}
\mathcal{S}(\gamma) = \sum_{\sigma,i} \langle \xi_i \rangle_\sigma^{P_\gamma} \xi_i^* \otimes q^{-c_1^{\text{vert}} (\sigma)} t^{-u_\gamma (\sigma)}
\end{equation}

Although we have defined the Seidel element differently in these two rings, note that the first determines the second, via the following lemma. 

\begin{lemma}\label{seidelequiv}
For any section class $\sigma \in P_\gamma$, there exists an additive homomorphism $\Phi_\sigma \co  \QH_* (X, \Lambda_\enr ) \to \QH_* (X, \Lambda_\univ)$ which takes $\mathcal{S}(\gamma, \sigma)$ to $\mathcal{S}(\gamma)$. This homomorphism restricts to the identity on $H_* (X)$.
\end{lemma} 

\noindent\textit{Proof} : Define 
\begin{equation*}
\Phi_\sigma (\xi_i \otimes e^{-A} ) = \xi_i \otimes q^{-c_1^{\text{vert}} (\sigma+A)}t^{-u_\gamma (\sigma+A)}.
\end{equation*}

Extend this map over $\QH_* (X, \Lambda_\enr )$ by linearity. This is clearly an additive homomorphism, and $\Phi_\sigma (\xi_i ) = \xi_i$. \proofend

We now explain what we mean when we say that the Seidel element is finite order.

\begin{definition}
Let $\Lambda$ be any Novikov ring, and let $\eta\in\QH_{2N}(X,\Lambda)$.  We say that $\eta$ has finite order if there exists $k$ such that $\eta^k = \mathbbm{1}\otimes\lambda$ for some $\lambda\in\Lambda$, $\lambda\neq 0$.
\end{definition}

This is not strictly the traditional sense of order, as some power is equal to $\mathbbm{1}\otimes\lambda$ rather than $\mathbbm{1}$.  However, by a result of Fukaya-Oh-Ohta-Ono \cite{fooo}, Lemma A.1, we know that any Novikov ring with coefficients in an algebraically closed field of characteristic $0$ is algebraically closed.  Therefore, by enlarging $\Lambda$ to have coefficients in $\C$ (we will call this $\Lambda^{\C}$), we can find $\eta\in\Lambda^{\C}$ such that $(\mathcal{S}(\gamma)\otimes\eta)^n = \mathbbm{1}$.  Therefore, the statement that $\mathcal{S}(\gamma)$ has finite order is true in the classic sense, up to multiplication by some $\eta\in\Lambda^{\C}$.  If $\mathcal{S}(\gamma)$ does not have finite order in this sense, we will say that it has \emph{infinite order}.

Note that the Seidel element $\mathcal{S}(\gamma)$ is in degree $2N$ for dimensional reasons.  But we can identify $\QH_{2N}(X,\Lambda_{\text{univ}}$) with the ring $\QH_{\ev}(X,\Lambda)$ by taking
\begin{equation*}
\psi(a\otimes q^{\epsilon_a}t^{\delta_a}) = a\otimes t^{\delta_a}
\end{equation*}

\noindent where $a\in H_{\text{ev}}(X)$, the subspace of $H_*(X)$ generated by even dimensional homology classes.  This map is an isomorphism, since $\psi^{-1}(a) = aq^{N - \frac{1}{2}\deg(a)}$.  Since $\Lambda$ is a field, working in this ring is more convenient for us.  Therefore, we will frequently use this isomorphism implicitly, especially in section \ref{nuvanishes} and Section \ref{grassmannians}.

Now we can define the action--Maslov homomorphism $I$ of Polterovich \cite{polterovich}. Although the original definition is the difference of the action functional and the Maslov class, Polterovich shows in \cite{polterovich}, Proposition 3.a, that the homomorphism can also be defined as the difference between the vertical Chern class and the coupling class. Namely,
\begin{equation}\label{actionmaslov}
u_\gamma = \kappa c_1^{\text{vert}} + I(\gamma) PD^{P_\gamma} (X). 
\end{equation}

Here, $\kappa$ is the same constant of monotonicity from before: $\omega = \kappa c_1(X)$.  We will use this alternate definition of the action--Maslov homomorphism, because it is more directly related to our results. Note in particular that if $\sigma$ is a section class with $c_1^{\text{vert}} (\sigma) = 0$, then $I(\gamma) PD^{P_\gamma} (X) = u_\gamma (\sigma)$. 

Finally, we will define Property $\mathcal{D}$. This should be seen as a generalization of the statement that the even degree homology classes of X are generated by divisors. Here, we will use the conventions that $\cdot$ represents the intersection product
\begin{equation*}
H_d(X)\otimes H_{2N-d}(X)\to H_0(X)\equiv\R,
\end{equation*}

\noindent and that all Gromov--Witten invariants are genus-zero invariants.  We will use these conventions throughout this paper.  Finally, we restate the definition of Property $\mathcal{D}$ from \cite{mcduff10}.

\begin{definition}\label{propertyd}(\cite{mcduff10})
$\QH_* (X, \Lambda)$ satisfies Property $\mathcal{D}$ if there exists an additive complement $\mathcal{V}$ in $H_{\text{ev}} (X, \Q)$ to the subring $\mathcal{D} \subset H_{\text{ev}} (X, \Q)$ generated by the divisors and the fundamental class such that: 
\begin{equation*}
\langle d, v\rangle_\beta=0 \text{ for all } \beta\in H_2^S(X)
\end{equation*}
\noindent for all $\mathcal{D} \ni d\neq [X]$, $v \in \mathcal{V}$.
\end{definition}


\section{Seidel Elements with Vanishing Valuation}\label{nuvanishes}
Let $(X, \omega)$ be a $2N$ dimensional symplectic manifold, $\gamma \in\pi_1 (\Ham(X, \omega))$, and $P_\gamma$ the bundle coming from the coupling construction. Quantum homology and Seidel elements in this section will always refer to those with respect to the universal Novikov ring $\Lambda$ defined in (\ref{lambdauniv}). 

In this section we will prove Proposition \ref{nuiszero} and thus Theorem \ref{ivanishes}. Proposition \ref{nuiszero} states that if $X$ has property $\mathcal{D}$ and every $\gamma$ has finite order $\mathcal{S}(\gamma)$, then $\nu(\lambda) = 0$. This in turn implies that $I=0$, which proves Theorem \ref{ivanishes}. 

The proof is adapted from the methods of McDuff \cite{mcduff10}.  McDuff proved a similar result - that the asymptotic spectral invariants descend - which is equivalent to $I=0$ for monotone symplectic manifolds.  Her assumptions were Property $D$ and a lower bound on the minimal Chern number, which in turn implies that $\mathcal{S}(n\gamma) = \mathbbm{1}\otimes\lambda + x$.  The shift in conditions requires some minor changes to the proof, which we present below.   When a lemma is cited as being from \cite{mcduff10}, this means that the same lemma was presented there, with conditions on the minimal Chern number instead of the Seidel element. If the proof is not given, then it proceeds exactly as in \cite{mcduff10}, using the lemmas from this paper in place of the originals.  The proofs that are presented use the same ideas as \cite{mcduff10}, but are modified in more significant ways.

We begin by defining a few specific terms which will we use throughout this section. 

\begin{definition} 
Let $Q_- (X, \Lambda) = \bigoplus_{i<2N} H_i (X) \otimes \Lambda$.
\end{definition}

\begin{definition}
Let $X$ and $\gamma$ be as above, and suppose that $\mathcal{S}(\gamma) = \mathbbm{1} \otimes \lambda + x$, where $x$ is any element in $Q_- (X, \Lambda)$ and $\lambda \neq 0$. Consider the sections $\lbrace\sigma\rbrace$ with $c_1^{\text{vert}} (\sigma) = 0$ which contribute to the Seidel element $\mathcal{S}(\gamma)$. Then define $\sigma_0$ to be a section such that $u_\gamma (\sigma_0 ) = \min\lbrace u_\gamma (\sigma)\rbrace$.
\end{definition}

\begin{note}
Since $\Lambda$ is a field, the condition that $\lambda \neq 0$ is equivalent to $\lambda$ being a unit in $\Lambda$.
\end{note}

The main thrust of our argument will be that knowing the Seidel element of $\gamma$ tells us a great deal about the Gromov--Witten invariants of $P_\gamma$. We use this knowledge to construct a homology representative of the Poincar\'{e} dual of $u_\gamma$, and to show that this homology representation has certain properties. 

\begin{lemma}\label{hlemma}
(Lemma 3.4, \cite{mcduff10}) Suppose $\mathcal{S}(\gamma) = \mathbbm{1} \otimes \lambda + x$, where $x$ is any element in $Q_- (X, \Lambda)$, and there is an element $H \in H_{2N} (P_\gamma )$ such that 
\begin{enumerate}
\item $H \cap [X]$ is Poincar\'{e} dual in $X$ to $[\omega]$. \label{hlemma1}
\item $H \cdot \sigma_0 = 0$. \label{hlemma2}
\item $H^{N +1} = 0$. \label{hlemma3}
\end{enumerate}
Then $\nu(\mathbbm{1} \otimes \lambda) = 0$.
\end{lemma}

 Here conditions (1) and (3) imply that $H$ is a representative of the Poincar\'{e} dual of the coupling class, so that (2) implies $\nu(\mathbbm{1} \otimes \lambda) = 0$. We wish to construct such an $H$. We will do so by ``fattening up'' a representative of the dual of $\omega$ in the fibre. As in \cite{mcduff10}, we define a map $s \co  H_* (X, \R) \to H_{*+2} (P, \R)$ by the identity 
 \begin{equation}
 s(a) \cdot_P v = \frac{1}{\langle pt\rangle^P_{\sigma_0}} \langle a,v\rangle^P_{\sigma_0},
 \end{equation}

\noindent for all $v \in H_* (P )$. Let $H_P = s(PD_X (\omega))$. Now we need to show that this $H_P$ satisfies the properties in Lemma \ref{hlemma}. Lemma \ref{fibervan} is a variant of parts (ii) and (iii) of \cite{mcduff10}, Lemma 4.2. Note that this version of the lemma eliminates the requirement on the minimal Chern number, but replaces it with a stronger condition on $\mathcal{S}(\gamma)$.  In turn, this gives us a stronger result; the lemma is true $\forall a\in H_*(X)$, rather than $a\in H_{<2N}(X)$ as in \cite{mcduff10}.

\begin{lemma}\label{fibervan} (based on Lemma 4.2 \cite{mcduff10})
 Suppose that $\mathcal{S}(\gamma) = \mathbbm{1} \otimes\lambda$ and let $\sigma = \sigma_0 -B$ for some $B \in H_2 (X)$ where $\omega(B) > 0$. Then for all $a \in H_* (X)$:
  \begin{enumerate}
  \item $\langle a, b\rangle^P_\sigma = 0, \forall b \in H_* (X)$.
  \item $\forall w \in H_* (P ), \langle a, w\rangle^P_\sigma$ depends only on $w \cap X$.
  \end{enumerate}
\end{lemma}

Despite these changes, the proof proceeds exactly as in \cite{mcduff10}, and will not be repeated here.

We will also need two lemmas from \cite{mcduff10}. The first is a special case of \cite{mcduff10}, Lemma 4.5. 

\begin{lemma}\label{sproperties}
Suppose that $\mathcal{S}(\gamma) = \mathbbm{1} \otimes \lambda$. Then: 
\begin{enumerate}
\item $s(pt) = \sigma_0$. 
\item $ s(a) \cap X = a, \forall a \in H_* (X,\R)$.
\end{enumerate}
\end{lemma}

The second is \cite{mcduff10}, Lemma 4.1. 

\begin{lemma} \label{sectioninvts}
Suppose that $a, b \in H_* (X), v, w \in H_* (P_\gamma,\R)$ and $B \in H_2 (X,\Z) \subset H_2 (P_\gamma,\Z )$. Then
\begin{enumerate}
\item$ \langle a, b, v\rangle^P_B = 0$.
\item $\langle a, v, w\rangle = \langle a, v \cap X, w \cap X\rangle^X_B$.
\end{enumerate}
\end{lemma}

For proofs of these two lemmas, see \cite{mcduff10}. Lemma \ref{sproperties} follows from simple computations using the definition of $s$:
\begin{equation*}
 s(pt) \cdot v = \sigma_0 \cdot v,
 \end{equation*}
 \noindent for any divisor class $v$. Lemma \ref{sectioninvts} essentially depends on the fact that the two fiber constraints can be located in different fibers. If $J$ is compatible with the fibration, the $J$-holomorphic curve will reside entirely in one fiber, and thus can intersect at most one of the fiber constraints. The second part follows similarly, since the $B$ curve must lie in the same fiber as $a$. 
 
The two hypotheses for the results in this section will be that $\QH_* (X,\Lambda)$ satisfies property $\mathcal{D}$ and that $\mathcal{S}(\gamma) = \mathbbm{1} \otimes \lambda$. We use Lemma \ref{divisors} at exactly two places in the proof of Proposition \ref{nuiszero} (specifically, in the proofs of Lemma \ref{twopointvan} and Lemma \ref{hpowervan}). The other results (and the main result, Proposition \ref{nuiszero}) thus require these conditions only so that they can use results of Lemma \ref{twopointvan} and Lemma \ref{hpowervan}. 

This lemma is very similar to $4.8.i$ in \cite{mcduff10}, but differs in a crucial way. Namely, the class $a$ here is permitted to be any class in $H_*(X)$, rather than being restricted to $H_{<2n}(X)$.  This case will be exactly the one needed in the proof of Lemma \ref{twopointvan}.  As such, the proof varies - the Gromov-Witten invariants vanish for different reasons - and is presented in its entirety.

\begin{lemma} \label{divisors}
Let $(X, \omega)$ be a symplectic manifold, $\gamma \in \pi_1 (\Ham(X, \omega))$, and $P$ the associated bundle. Assume $\mathcal{D}$ is the part of $H_{\text{ev}} (P,\R )$ generated by divisors $\lbrace D_i , X\rbrace$ and $\mathcal{S}(\gamma) = \mathbbm{1} \otimes \lambda$. Let $a \in H_* (X,\R)$ be a fiber class, $v \in \mathcal{D}$, and $B \in H_2 (X,\Z)$ such that $\omega(B) > 0$. Then the Gromov--Witten invariant $\langle v, a \rangle^P_{\sigma_0 -B}$ vanishes.
\end{lemma}

\noindent \textit{Proof} : Suppose not. As in the proof of 4.8.i of \cite{mcduff10}, we take a section of minimal energy such that some invariant of this form does not vanish and call it $\sigma' $. $v$ is a product of divisors, and we claim we may assume that each of these divisors $D_i$ satisfies $D_i \cdot \sigma ' = 0$. First, we can show that none of the $D_i = X$.  If any of them did, then $v$ would be a fiber class, and $\langle v, a \rangle^P_{\sigma_0 - B}$ would vanish by Lemma \ref{fibervan} (note that we need here the stronger condition of Lemma \ref{fibervan}, since we could have $a\in H_{2N}(X)$).  Then to any other $D_i$,  we can add a multiple of $X$ to obtain a new class $D_i' $ which differs from $D_i$ by a fiber class and has $D_i' \cdot \sigma ' = D_i \cdot \sigma' + kX \cdot \sigma' = 0$, for appropriate choice of $k$. Lemma \ref{fibervan} shows that adding a fiber class to v does not change our Gromov--Witten invariant. Now consider the set 
\begin{equation*}
\lbrace v_i | \langle v_i , a\rangle^P_{\sigma'} \neq 0\rbrace.
\end{equation*}

Each of these $v_i$ is a linear combination of products of $k$ divisors.  We will assume, without loss of generality, that  v is one of these $v_i$ and it is exactly a product of $k$ divisors.   We will perform induction on $k$.

If $k = 1$, the invariant vanishes by the divisor axiom for Gromov--Witten invariants (see \cite{ms}, Section 7.5), which says that $\langle D_1 , a\rangle^P_{\sigma'} = \langle a  \rangle^P_{\sigma'} (D_1\cdot \sigma ' ) = 0$. If $k > 1$, we use Theorem 1 of Lee--Pandharipande from \cite{lp}, as restated in \cite{mcduff10}, (4.2). This identity is stated as follows. Take a basis $\xi_i$ of $H_* (X)$ and extend it to a basis of $H_* (P)$ by adding classes $\xi_i^*$ such that $\xi_i \cdot \xi^*_j = \delta_{ij}$ and $\xi_i^* \cdot \xi_j^* = 0$ (note that these $\xi_i^*$ are not the same as the $\xi_i^*$ above, as they form a dual basis in $H_*(P)$ rather than $H_*(X)$). Note that the $\xi_i$ here are fiber classes, but the $\xi_i^*$ cannot be fiber classes. Now take classes $u, v, w \in H_* (P_\gamma ), H \in H_{2N} (P )$ a divisor, and $\alpha \in H_2 (P )$. Then Lee and Pandharipande show that 
\begin{multline}\label{leepanidentity}
\langle Hu, v,w\rangle^P_\alpha = \langle u, Hv,w\rangle^P_\alpha + (\alpha \cdot H) \langle u, \tau v,w\rangle^P_\alpha\\
 - \sum_{i,\alpha_1 +\alpha_2 =\alpha}(\alpha_1 \cdot H) \langle u, \xi_i,\ldots\rangle^P_{\alpha_1} \langle \xi_i^*, v,\ldots\rangle^P_{\alpha_2}
\end{multline}

\noindent where $\tau$ is a descendant constraint and ``$\ldots$'' indicates that the $w$ term may appear in either factor. 

Now, assume that the statement is true for all $v\in\mathcal{D}$ of codimension $2k-2$. Let $v = D_1 \cdots D_{k-1} \cdot D_k$ (where, as above, we can assume that all of these divisors have $D_i \cdot \sigma ' = 0$). Given any section class $\sigma$, McDuff shows in Lemma 2.9 of \cite{mcduff00} that in the above sum, a section class can only decompose into $\sigma -\alpha$ and $\alpha$ where either $\alpha$ is a fiber class or $\sigma -\alpha$ is a fiber class. In both cases, the other element of the decomposition will be a section class by necessity. This follows from considering $J$-holomorphic curves where $J$ is compatible with the fibration. By combining equation \ref{leepanidentity} and this decomposition into divisors, one sees that (we take $w = D_1 \cdots D_{k-1}$ and $D = D_k$, to simplify our notation) 
\begin{align}
\langle w \cdot D, a\rangle^P_{\sigma'} &= \langle w \cdot D, a,X\rangle^P_{\sigma'}\\
&= \langle w, D \cdot a,X\rangle^P_{\sigma'} + (D \cdot \sigma ' ) \langle w, \tau a,X\rangle^P_{\sigma'}\label{lp1}\\
&-\sum_{c\in H_2 (X),i}((\sigma' - c)\cdot D)\langle w,\xi_i,\ldots\rangle^P_{\sigma' - c}\langle \xi_i^*,a,\ldots\rangle^P_c\label{lp2}\\
&-\sum_{c\in H_2 (X),i}((\sigma' - c)\cdot D)\langle w,\xi_i^*,\ldots\rangle^P_{\sigma' - c}\langle \xi_i,a,\ldots\rangle^P_c\label{lp3}\\
&-\sum_{c\in H_2 (X),i}(c\cdot D)\langle w,\xi_i,\ldots\rangle^P_{c}\langle \xi_i^*,a,\ldots\rangle^P_{\sigma' - c}\label{lp4}\\
&-\sum_{c\in H_2 (X),i}(c\cdot D)\langle w,\xi_i^*,\ldots\rangle^P_{c}\langle \xi_i,a,\ldots\rangle^P_{\sigma' - c}\label{lp5}.
\end{align}

We will go through the right hand side of this equation line by line and show that each of them must vanish. Line (\ref{lp1}) has two terms -- the first one vanishes because $w$ is of codimension $2k-2$ and the second one vanishes because $D \cdot \sigma ' = 0$. If line (\ref{lp2}) does not vanish then we must have either $c = 0$, or $\omega(c) > 0$ and $X$ in the first factor (otherwise the second factor would vanish by Lemma \ref{sectioninvts}). If $c = 0$, then $(\sigma ' - c) \cdot D = \sigma ' \cdot D = 0$ and line (\ref{lp2}) vanishes. Thus our first factor is 
\begin{equation*}
\langle w, \xi_i,X\rangle^P_{\sigma '-c}
\end{equation*}

\noindent with $\omega(c) > 0$, which vanishes by the minimality of $\sigma ' $. Line (\ref{lp3}) must vanish by Lemma \ref{sectioninvts} because the second factor is a fiber invariant with two fiber constraints. In line (\ref{lp4}), the $X$ must insert into the second term by Lemma \ref{sectioninvts}, and thus we have invariants of the form 
\begin{equation*}
\langle w, \xi_i\rangle^P_c \langle \xi_i^* , a,X\rangle^P_{\sigma '-c}.
\end{equation*}

Note that $c \neq 0$, so this vanishes by minimality of $\sigma ' $.

Finally, line (\ref{lp5}) must vanish because the second factor is of the form $\langle \xi_i, a,X\rangle^P_{\sigma ' -c} = \langle \xi_i, a \rangle^P_{\sigma ' -c}$. This invariant vanishes because it has two fiber constraints.  Assume that it does not vanish.  Then it would contribute to $\mathcal{S}(\gamma) * a$, as in\cite{ms}, (11.4.4), 
\begin{equation}
\mathcal{S}(\gamma) * a = \sum_{i,\sigma} \langle a, \xi_i\rangle^{P_\gamma}_\sigma \xi_i^*\otimes q^{-c_1^{\text{vert}} (\sigma)}t^{-u_\gamma (\sigma) }.
\end{equation}

But since it doesn't vanish, $\omega(c) > 0$, and thus $\sigma ' - c$ has less energy than $\sigma_0$, which contradicts the definition of $\sigma_0$. Therefore, the entire invariant vanishes, and by induction, all such invariants vanish.\proofend

The following lemma differs from Lemma 4.6 in \cite{mcduff10} only in the initial conditions, and the proof follows in the same way, using the modified lemmas where appropriate.  The one significant difference from \cite{mcduff10} is highlighted.
\begin{lemma} \label{twopointvan}
Assume the conditions of Proposition \ref{nuiszero}. Then $\langle h, \sigma_0\rangle^P_{\sigma_0} = 0$ where $h \in H_* (X,\R)$ is the Poincar\'{e} dual of the symplectic form in $X$. 
\end{lemma}

\noindent\textit{Proof} : We can take any divisor class $D$ in $P$ such that $D \cap X = h$ and add copies of $X$ to get a class $K$ such that $K \cap X = h$ and $K \cdot \sigma_0 = 0$. Then the identity of Lee--Pandharipande gives us 
\begin{align*}
\langle h, \sigma_0\rangle^P_{\sigma_0} &=  \langle h,\sigma_0,X\rangle^P_{\sigma_0}\\
&=\langle X, K\sigma_0,X\rangle^P_{\sigma_0} + (\sigma_0 \cdot K) \langle X, \tau \sigma_0,X\rangle^P_{\sigma_0}\\
&- \sum_{\alpha\in H_2 (P )} ((\sigma_0 - \alpha) \cdot K)\langle X,\xi_i,\ldots\rangle^P_{\sigma_0 - \alpha}\langle \xi_i^*,\sigma_0,\ldots\rangle^P_\alpha\\
&- \sum_{\alpha\in H_2 (P )} ((\sigma_0 - \alpha) \cdot K) \langle X,\xi_i^*,\ldots\rangle^P_{\sigma_0 - \alpha}\langle \xi_i,\sigma_0,\ldots\rangle^P_\alpha
\end{align*}

The proof of how nearly all of these invariants vanish proceeds as in \cite{mcduff10}, and will not be repeated.  The one case which differs is that of the invariants in the second sum when $\sigma_0 - \alpha$ is a section class. The invariant in question is thus of the form: 
\begin{equation*}
\langle X, \xi_i^*,X\rangle^P_{\sigma_0-\alpha}\langle\xi_i , pt\rangle^X_\alpha
\end{equation*}

\noindent with $\alpha \neq 0$ and $\omega(\alpha) > 0$ since $((\sigma-\alpha)\cdot K) = 0$. If the second factor does not vanish, then property $\mathcal{D}$ tells us that the class $\xi_i$ is generated by divisors, because $pt \in \mathcal{D}$. If $\xi_i$ is generated by divisors, $\xi_i^*$ must be generated by divisors, and Lemma \ref{divisors} tells us that the first factor must vanish.  Note that here that Lemma $4.8.i$ from \cite{mcduff10} would not be sufficient, since the fiber class here is $X$.\proofend

\begin{corollary}\label{hdotsigma0} (Corollary 4.7, \cite{mcduff10})
Assuming the conditions of Proposition \ref{nuiszero}, we have $H \cdot \sigma_0 = 0$.
\end{corollary}

\begin{lemma} \label{hpowervan} (Lemma 4.8, \cite{mcduff10})
Assuming the conditions of Theorem \ref{ivanishes}, 
\begin{equation*}
\langle H^{N +1-k} , X\cap H^k\rangle^P_{\sigma_0} = 0
\end{equation*}
for all $k$. 
\end{lemma}

The proof of this statement follows exactly as Lemma 4.8 in \cite{mcduff10}, and the following corollary is identical to Corollary 4.9 in \cite{mcduff10}. 

\begin{corollary} \label{hnplus1van} (Corollary 4.9, \cite{mcduff10})
Assuming the conditions of Proposition \ref{nuiszero}, $H^{N +1} = 0$. 
\end{corollary}

The proof of Proposition \ref{nuiszero} now follows from these results.

\noindent\textit{Proof of Proposition \ref{nuiszero}}: $X$ satisfies the conditions of Lemma \ref{twopointvan}. Then by Corollary \ref{hdotsigma0}, $H_{P_\gamma} \cdot \sigma_0 = 0$. Similarly, Lemma \ref{hnplus1van} shows that $H^{m+n+1} = 0$. Thus, $H$ satisfies the conditions of the Lemma \ref{hlemma}. \proofend

\noindent\textit{Proof of Theorem \ref{ivanishes}}: Let $[\gamma]\in\pi_1(\Ham(X,\omega))$. Since $S(n\gamma)$ is of the form $\mathbbm{1} \otimes \lambda$ and $\QH_*(X,\Lambda)$ has Property $\mathcal{D}$, Proposition \ref{nuiszero} implies that $\nu(\mathcal{S}(n\gamma)) = 0$. This implies that there exists a section class $\sigma_0$ which contributes to $\mathcal{S}(n\gamma)$ and has $u_{n\gamma}(\sigma_0) = 0$.  By the definition of the Seidel element, $c_1^\text{vert}(\sigma_0)$ is also $0$, and thus $c_1^{\text{vert}} (\sigma_0 ) = u_{n\gamma} (\sigma_0 ) = 0$. Therefore, by Equation \ref{actionmaslov}, $I(n \gamma  ) = 0$. But $I(n \gamma ) = nI(\gamma)$, so $I(\gamma) = 0$. \proofend

\section{Manifolds with $I=0$}\label{examples}
Now we will discuss several monotone symplectic manifolds that we can show have $I = 0$. We do this by showing results about property $\mathcal{D}$ and the form of the Seidel element. Note that we need to show that the Seidel element $S(k\gamma ) = \mathbbm{1} \otimes \lambda$. In some cases, it may be easier to show that the enriched Seidel element $S(k\gamma , k\sigma)$ has this form. Then Lemma \ref{seidelequiv} implies that $S(k\gamma )$ has the form $\mathbbm{1} \otimes \lambda$.

We begin with a lemma which helps us to reduce the potential elements of quantum homology which can be Seidel elements.
\begin{lemma}
If $(X, \omega)$ is monotone and $\gamma \in \pi_1 (\Ham(X, \omega))$, then the Seidel element $\mathcal{S}(\gamma, \sigma')$ (and also the Seidel element $\mathcal{S}(\gamma))$ will have coefficients in $\Z$ and a finite number of terms. 
\end{lemma}

This result is very straightforward and uses very standard techniques in Gromov-Witten theory.  A full proof can be found in \cite{branson}.  Using this lemma, we need only consider elements in a smaller subring of $QH_*(X,\Lambda_\enr)$:

\begin{definition}
Let $\Lambda_{\enr,\Z}$ be the subring of $\Lambda_\enr$ with integral coefficients. Then define $Q_\enr (X) = \QH_* (X, \Lambda_{\enr,\Z} )$ to be the subring of $\QH_* (X, \Lambda_\enr )$ which consists of finite sums of elements with coefficients in $\Z$. Thus, a typical element is
\begin{equation*}
\sum^n_{i=0} x_i e^{-C_i}
\end{equation*}

\noindent where $x_i \in H_* (X, \Z)$ and $C_i \in H_2^S (X, \Z)$.
\end{definition}

\subsection{$\CP^m \times \CP^n$}\label{cpnsection}
Let $X$ be $\CP^m \times \CP^n$ with the monotone symplectic form, and let $N = m + n$.

\begin{note}
In fact, all of the results in this section are true for products of an arbitrary number of projective  spaces with the monotone symplectic form.  For simplicity of notation, though, we will prove them for $\CP^m\times\CP^n$ only.
\end{note}

In order to show that $I = 0$, we need to show that some power of the Seidel element is of the form $\mathbbm{1} \otimes \lambda$. This is a consequence of the algebraic structure of the quantum homology: namely, that the subring $Q_\enr$ is an integral group ring over an ordered group. 

\begin{definition}
 An ordered group is a group $G$ equipped with a total order $<$ which is translation invariant: $g < h \Rightarrow g \cdot a < h \cdot a$ and $a \cdot g < a \cdot h$ $\forall g, h, a \in G$. 
 \end{definition}
 
 \begin{theorem}\label{cpnprodordered}
If $X$ is $\CP^m\times\CP^n$ with the monotone symplectic form, then $Q_\enr (X)$ is an $\Z$ group ring over an ordered group.
 \end{theorem}

Theorem \ref{cpnprodordered} follows directly from Lemmas \ref{prodordered} and \ref{cpnordered}. 

\begin{lemma}\label{prodordered}
If $Q_\enr (X ' , \Lambda_{\enr,\Z}^{X'} )$ and $Q_\enr (X '' , \Lambda_{\enr,\Z}^{X''})$ are both integral group rings over ordered groups, then $Q_\enr (X ' \times X '' , \Lambda_{\enr,\Z}^{X'\times X''} )$ is also an integral group ring over an ordered group.
\end{lemma}

\noindent\textit{Proof} : First, note that $H_2^S(X'\times X'') \cong H_2^S(X')\oplus H_2^S(X'')$.  Therefore, we can write any $e^A$ for $A\in H_2^S(X'\times X'')$ as $e^{A'}e^{A''}$ where $A'\in H_2^S(X')$ and $A''\in H_2^S(X'')$.  This gives us an isomorphism 

\begin{equation*}
\Lambda_{\enr,\Z}^{X'\times X''} \cong \Lambda_{\enr,\Z}^{X'}\otimes\Lambda_{\enr,\Z}^{X''}
\end{equation*}

Therefore, we have $Q_\enr (X ' \times X '' , \Lambda_{\enr,\Z}^{X'\times X''} )\cong Q_\enr (X ' \times X '' , \Lambda_{\enr,\Z}^{X'} \otimes \Lambda_{\enr,\Z}^{X''} ) \cong H_* (X ' \times X '' ) \otimes_\Z \Lambda_{\enr,\Z}^{X'} \otimes \Lambda_{\enr,\Z}^{X''} $ as additive groups. Thus, by the classical K\"{u}nneth formula, $Q_\enr (X ' \times X '' , \Lambda_{\enr,\Z}^{X'} \otimes \Lambda_{\enr,\Z}^{X''} ) \cong Q_\enr (X ' , \Lambda_{\enr,\Z}^{X'} ) \otimes_\Z QH_\enr (X '' , \Lambda_{\enr,\Z}^{X''} )$ as additive groups. By the quantum Kunneth formula (Section 11.1, \cite{ms}),  this is actually a ring isomorphism. Therefore, $Q_\enr (X ' \times X '' , \Lambda_{\enr,\Z}^{X'} \otimes \Lambda_{\enr,\Z}^{X''} ) \cong Q_\enr (X ' , \Lambda_{\enr,\Z}^{X'} ) \otimes_\Z Q_\enr (X '' , \Lambda_{\enr,\Z}^{X''} )$.

We assumed that both of these subrings were integral group rings over ordered groups. Therefore, we have that $Q_\enr (X ' \times X '' , \Lambda_{\enr,\Z}^{X'} \otimes \Lambda_{\enr,\Z}^{X''} ) \cong \Z(G' ) \otimes_\Z \Z(G'' )$. But this is isomorphic to $\Z(G' \times G'' )$. Give $G' \times G''$ the lexicographic ordering. The product of two ordered groups with the lexicographic ordering is still an ordered group: $(g ', g '' ) < (h' , h'' ) \Rightarrow (g ' \cdot a' , g '' \cdot a'' ) < (h' \cdot a' , h'' \cdot a'' )$ and $(a' \cdot g ', a'' \cdot g '' ) < (a' \cdot h' , a'' \cdot h'' )$. \proofend

Note that this statement is NOT true with the universal coefficients $\Lambda_{\text{univ}}$.  In that case, the tensor product is over $\Z[q]$ and the isomorphism does not respect the ordering.  Therefore, the enriched coefficients are necessary here to obtain the desired result.

\begin{lemma}\label{cpnordered}
$Q_\enr (\CP^n )$ is a group ring over an ordered group.
\end{lemma}

\noindent\textit{Proof}: 
Recall that 
\begin{equation*}
Q_\enr (\CP^n ) \cong \frac{\Lambda_{\enr,\Z} [x]}{\langle x^{n+1} = e^{-A}\rangle},
\end{equation*}

\noindent where $A$ is the class of the generator in $H_2(X,\Z)$.  Let $q = e^A$ and then let $G$ be the group generated by $x$ and $q$ with relation $x^{n+1} = q^{-1}$ (note that $(G,*)$ is isomorphic to $(\Z,+)$, via the isomorphism $\phi(x^k) = k$). This group can be ordered by using the mapping $\phi \co  (G, \cdot) \to (\Q, +)$ where $\phi(x) = \frac{1} {n+1}$ and $\phi(q) = -1$. Then $G$ is ordered by the pullback of the ordering on $\Q$. Clearly, $Q_\enr (\CP^n )$ is just the $\Z$ group ring of $G$.\proofend

Now we can combine Theorem \ref{cpnprodordered} and an algebraic lemma to determine the units of $Q_\enr (X)$.

\begin{lemma}\label{orderedunits}
If $G$ is an ordered group, then the units of $\Z(G)$ are $\pm G$.
\end{lemma}

Lemma \ref{orderedunits} is proved as Lemma 45.3 of \cite{sehgal}. The proof in Sehgal is incomplete, so we provide a corrected version here: 

\noindent\textit{Proof (Sehgal)}: Take a nonmonomial unit $p = \sum_{i=1}^t u_i * g_i$ of the group ring and its inverse (which must also then be nonmonomial) $p^{-1} = \sum_{i=1}^\ell v_i * h_i$, with $g_1 < g_2 < \ldots < g_t$ and $h_1 < h_2 < \ldots < h_\ell$. If we multiply these two elements, we get $1_G = u_1 v_1 * g_1 h_1 + \ldots + u_t v_\ell * g_t h_\ell$. Then, for this equation to be true, the group element in any term on the right hand side must be $1_G$ or cancel with the group element from another term. Since $G$ is ordered, $g_1 h_1 < g_i h_j$, for $i > 1$ or $j > 1$ and $g_i h_j < g_t h_\ell$ for $i < t$ or $j < \ell$, so these group elements cannot cancel with other terms. Thus, we must have $g_1 h_1 = 1_G = g_t h_\ell$ and thus $g_1^{-1} = h_1$ and $g_t^{-1} = h_\ell$. But we have $g_1 < g_t \Rightarrow g_1^{-1} > g_t^{-1} \Rightarrow h_1 > h_\ell$, which is a contradiction. Therefore, $p$ must be monomial.\proofend

\begin{corollary} \label{cpnprodmon}
Let $X = \CP^m \times \CP^n$ with the monotone symplectic form. Then the only units in $Q_\enr (X)$ are the monomial units - those of the form $\pm a^i b^j \otimes e^C$, $C \in H_2 (X, \Z)$.
\end{corollary}
\noindent\textit{Proof}: Theorem \ref{cpnprodordered} shows that $Q_\enr (X)$ is isomorphic to an integral group ring over the group generated by $a$, $b$, and $e^A$ (where $A$ is a generator of $H_2 (X, \Z)$). Since this group is ordered, all of its units are monomial by Lemma \ref{orderedunits}. \proofend

\begin{theorem} \label{cpnprodcyclic}
For $\CP^m \times\CP^n$ with the monotone symplectic form and for any loop $\gamma \in \pi_1 (\Ham(\CP^m \times \CP^n , \omega))$, $\mathcal{S}(\gamma)$ has finite order.
\end{theorem}

\noindent\textit{Proof} : Let $\sigma$ be a section class in $H_2 (P_\gamma )$. Corollary \ref{cpnprodmon} shows that $\mathcal{S}(\gamma, \sigma)$ must be of the form $\pm a^f b^g \otimes e^{-C}$. Let $k = (m + 1)(n + 1)$. Then $S(k\gamma , k\sigma) = a^{kf} b^{kg} \otimes e^{-C} = \mathbbm{1} \otimes e^{(kh-(n+1)f )A} e^{-C}$. By Lemma \ref{seidelequiv}, the same $k$ also works for $\mathcal{S}(\gamma k ) = \mathbbm{1} \otimes \lambda$. \proofend

\noindent\textit{Proof of Theorem \ref{cpnproducts}}: The first condition of Theorem \ref{ivanishes} is satisfied because all classes in $H_* (\CP^{n_1} \times \ldots \times \CP^{n_k} )$ are generated by divisors. Therefore, $\QH_* (X, \Lambda)$ satisfies property $\mathcal{D}$. Theorem \ref{cpnprodcyclic} shows that the second condition is satisfied for a product of two projective spaces, and thus Theorem \ref{ivanishes} shows that $I = 0$. By using Lemma \ref{prodordered} $(k - 1)$ times, one can show Theorem \ref{cpnprodcyclic} for the product of $k$ projective spaces. The result follows. \proofend

\subsection{$G(2, 4)$}\label{grassmannians}

The complex grassmannians are another class of monotone symplectic manifolds with well-understood quantum homology.  Let $(X,\omega)$ be the Grassmannian of $2$-planes in $\C^4$ (which we will also denote by $G(2,4)$) with $\omega = c_1(X)$.  Theorem \ref{g24} states that the action--Maslov homomorphism vanishes for $(X,\omega)$, the simplest grassmannian which is not a projective space. We need to show two things to prove this statement: that $\mathcal{S}(k\gamma) = \mathbbm{1} \otimes \lambda$ and that $X$ satisfies property $\mathcal{D}$. First, we will show that the Seidel element must have finite order. Unlike the products of projective spaces, here we do not even need to use enriched coefficients. Instead of $Q_\enr (X)$, we will look at the analogous subring $Q(X)$ of finite sums with integral coefficients in $\QH_* (X, \Lambda)$.  Let $x_1$ and $x_2$ be the Poincar\'{e} duals of the first and second Chern classes, respectively.  Then Siebert and Tian \cite{st97} show that the ring $Q(X)$ is 
\begin{equation*}
\frac{\Lambda[x_1 , x_2 ]}{\langle x_1^3 - 2x_1 x_2 , x_1^2 x_2 - x_2^2 - t^{4}\rangle }.
\end{equation*}

Because $\dim(X) = 8$ and the minimal Chern number is $4$, the terms which can appear in the Seidel element are sharply limited. Assume that a section $\sigma$ contributes to the Seidel element. Then any other contributing sections are of the form $\sigma' = \sigma + kL$, where $k \in \Z$ and $L = x_1 x_2$ is the class of a line in $X$. Since contributing sections must have $-8 \leq c_1^{\text{vert}} (\sigma ' ) \leq 0$, clearly another section can only exist if $\sigma ' = \sigma \pm L$ and $c_1^{{\text{vert}}} (\sigma) = 0$ or $-8$. $H_* (X)$ has generators organized by degree as follows: 
\begin{center}
\begin{tabular}{c c c c c}
$0$ & $2$ & $4$ & $6$ & $8$ \\
$x_2^2$ & $x_1x_2$ & $x_1^2, x_2$ & $x_1$ & $\mathbbm{1}$
\end{tabular}
\end{center}

The Seidel elements form a subgroup of the units - the product of two Seidel elements is a Seidel element, and so is the inverse. All of these elements have degree equal to the dimension of $X$, which is $8$. Thus, the Seidel element can only be of the form: 
\begin{align}
&a\mathbbm{1} t^\epsilon + bx_2^2 q^4 t^{4+\epsilon}\label{firsttype24}\\
&a x_1 q^1 t^\epsilon \label{secondtype24}\\
&ax_1^2 q^2 t^\epsilon + bx_2 q^2 t^\epsilon\label{thirdtype24}\\
&ax_1 x_2 q^3 t^\epsilon\label{fourthtype24}
\end{align}

Since these elements are of degree $8$ in $\QH_*(X,\Lambda_\univ)$, we will work with coefficients in $\Lambda$ instead. Similarly, the exponent of $t$ is determined up to a constant multiple $\lambda = t^\epsilon$ so we will also suppress $t$. These elements must be units in the quantum homology (with inverses of the same form), and since the symplectic form is monotone, $a$ and $b$ must be integers.   

\begin{lemma}\label{g24units}
The Seidel element, up to appropriate powers of $q$ and $t$, is either $\pm \mathbbm{1}, \pm x_2 , \pm (x_1^2 - x_2)$, or $\pm x_2^2$.
\end{lemma}

\noindent\textit{Proof} : We will proceed by checking each possible Seidel element, starting with (\ref{secondtype24}).  In this case, $\mathcal{S}(\gamma^{-1} )$ is of the form given in (\ref{fourthtype24}). Therefore, with appropriate powers of $q$, we have
\begin{align*}
\mathbbm{1} &= (ax_1 q^{-1} ) * (bx_1 x_2 q^{-3})\\
&= abx_1^2 x_2 q^{-4}\\
& = ab(\mathbbm{1} + x_2^2 q^{-4}).
\end{align*}

This implies that $ab = 0$ and $ab = 1$, which is impossible. Thus, no such elements can be Seidel elements. Now we look at (\ref{firsttype24}). First, note that
\begin{align*}
x_2^4 &= (x_2^2)^2\\
	   &= (q^4 - x_1^2 x_2)^2\\
	   &= q^8 - 2x_1^2x_2q^4 + x_1^4x_2^2\\
	   &= q^8 - 2x_1^2x_2(x_1^2x_2 - x_2^2) + x_1^4x_2^2\\
	   &= q^8 - x_1^4x_2^2 + 2x_1^2 x_2^3\\
	   &= q^8 + x_1x_2^2(2x_1x_2 - x_1^3)\\
	   &= q^8
\end{align*}

Therefore,  we have
\begin{align*}
\mathbbm{1} &= (a\mathbbm{1} + bx_2^2q^{-2} ) * (c\mathbbm{1} + dx_2^2 q^{-2}) \\
	&= (ac + bd)\mathbbm{1} + (ad + bc)x_2^2.
\end{align*}

Hence, $(ad + bc) = 0$ and $(ac + bd) = 1$. Then, either $d = 0$ or $a =\frac{-bc}{d}$.  First, we will address the case where $d\neq 0$.  By substituting $\frac{-bc}{d}$ for $a$, one obtains that $b = \frac{d}{ d^2 -c^2}$ and $a = \frac{-c} { d^2 -c^2}$. Since $a$ and $b$ are both integers, this means that $d^2 - c^2$ divides both $c$ and $d$. Since $d\neq 0$, this is only true when $\lbrace a, b, c, d\rbrace = \lbrace 0, \pm1, 0, \pm1\rbrace$.  If $d=0$, the equations immediately show that we must have $\lbrace\pm1, 0, \pm1, 0\rbrace$. Therefore, the Seidel element must be either $\mathbbm{1}$ or $x_2^2$ multiplied by some $\lambda$.  

Finally, we look at (\ref{thirdtype24}). Here, we will have
\begin{align*}
\mathbbm{1} &= (ax_1^2q^{-2} + bx_2q^{-2}  ) * (cx_1^2q^{-2}  + dx_2q^{-2}  ) \\
	&= acx_1^4q^{-4}  + (bc + ad)x_1^2 x_2q^{-4} + bdx_2^2q^{-4}\\
	&= (ac+bc+ad)*(2x_1^2x_2q^{-4}) + bd x_2^2q^{-4}\\
	&= (ac + bc + ad)*(2x_2^2q^{-4} + 2*\mathbbm{1}) + bd x_2^2q^{-4}\\
	& = (2ac + bc + ad + bd)x_2^2q^{-4} + (2ac + bc + ad)\mathbbm{1}.
\end{align*}

This will be true if and only if $(2ac + bc + ad) = -bd = 1$. Thus we have $b = -d = \pm1$. If $b = -d = 1$, then $c = \frac{1+a}{1+2a}$, which is only integral if $\lbrace a, b, c, d\rbrace = \lbrace 0, 1, 1, -1\rbrace$ or $\lbrace -1, 1, 0, -1\rbrace$. Similarly, if $b = -d = -1$, then $c= \frac{1-a}{-1+2a}$, which is only integral if $\lbrace a, b, c, d\rbrace = \lbrace 0, -1, -1, 1\rbrace$ or $\lbrace 1, -1, 0, 1\rbrace$. Therefore, we have either $\pm(x_1^2 - x_2 )$ or $\pm x_2$. This completes the proof. \proofend

\begin{lemma} \label{g24cyclic}
Let $\mathcal{S}(\gamma)$ be the Seidel element of $\gamma \in \pi_1 (\Ham(X, \omega))$. Then $S(4\gamma ) = \mathbbm{1} \otimes \lambda$.
\end{lemma}

\noindent\textit{Proof}: $\mathcal{S}(\gamma)$ must be of a form listed in Lemma \ref{g24units}. Clearly, since $\Lambda$ is a field, the coefficient $\lambda$ does not affect invertibility, and we only need to concern ourselves with the homology terms. Since $\mathbbm{1}^4 = \mathbbm{1}$ and $x_2^4 = q^8$, this is obvious for the first two cases. This leaves only the case where $\mathcal{S}(\gamma) = x_1^2 -x_2$. But $(x_1^2 - x_2 )^2 = x_2^2$ (by the calculations for line (\ref{thirdtype24})), so the statement also holds in this case. \proofend

In order to show that the action--Maslov homomorphism vanishes on $G(2, 4)$, we also need to show that it satisfies Property $\mathcal{D}$. This is slightly weaker than requiring that the quantum homology be generated by divisors. $\QH_*(G(2, 4), \Lambda)$ is not generated by divisors, but does satisfy property $\mathcal{D}$. 

\begin{lemma} \label{g24d}
The quantum homology of $G(2, 4)$ with coefficients in $\Lambda$ satisfies property $\mathcal{D}$. 
\end{lemma}
\noindent\textit{Proof}: First, note that the homology is generated (over $\Q$) by $x_1$ in every degree except $4$ ($x_1 x_2 = \frac{1}{2} x_1^3$ and $x_2^2 = x_1^2 x_2 = \frac{1}{2} x_1^4$ ). Therefore, $\mathcal{V}$ must be generated by some class $ax_1^2 + bx_2$. But by the dimension formula for genus $0$ Gromov--Witten invariants, if 
\begin{equation*}
\langle d, v\rangle_A\neq 0,
\end{equation*}

\noindent then the codegrees of $d$ and $v$ must add up to $8 + 2c_1 (A) - 2$. If $A \neq 0$, we have $c_1 (A) \geq 4$, so the sum of the codegrees must be at least $14$. But if $v$ is $ax_1^2 + bx_2$, it has codegree $4$ and $d$ must have codegree $10$. But $G(2, 4)$ is $8$ dimensional, so this cannot happen. Therefore, $\langle d, ax_1^2 + bx_2\rangle_A = 0$ for all $d\in\mathcal{D}$ and $G(2, 4)$ satisfies property $\mathcal{D}$. \proofend

Lemmata \ref{g24cyclic} and \ref{g24d} are sufficient to show that the action--Maslov homomorphism vanishes on $G(2, 4)$.  This completes the proof of Theorem \ref{g24}.

\subsection{Other Grassmannians}
The immediate question posed by Theorem \ref{g24} is whether these results can be extended to other grassmannians.  The answer to this question is unfortunately no.  We can show, using the quantum Schubert calculus developed by Bertram \cite{bertram97}, that $G(2,4)$ and $\CP^n$ are the only grassmannians with property $\mathcal{D}$.  Additionally, we can show that for $G(2,2n+1)$, there exist units in quantum homology with infinite order.  Both of these results require an extensive treatment of the quantum Schubert calculus - see \cite{branson} for more details.  Because of these results, Theorem \ref{nuiszero} will not suffice to show that $I=0$ for these grassmannians.

However, it is possible to replace property $\mathcal{D}$ with a weaker statement about the quantum homology in Theorem \ref{nuiszero}, as property $\mathcal{D}$ requires the vanishing of far more Gromov--Witten invariants than is actually necessary.  It is possible that doing so could provide progress on the question of whether $I$ vanishes on the higher dimensional grassmannians.

\bibliographystyle{amsalpha}
\bibliography{actionmaslov_final}

\end{document}